\documentclass{amsart}

\usepackage{amssymb}
\usepackage{amsmath}
\usepackage{amsthm}
\usepackage{graphicx,psfrag}

\theoremstyle{plain}
\newtheorem{theorem}{Theorem}
\newtheorem{lemma}[theorem]{Lemma}

\theoremstyle{definition}
\newtheorem{definition}[theorem]{Definition}
\newtheorem*{notation}{Notation}

\theoremstyle{remark}

\newcommand{\Int}{\mathrm{int}}

\begin{document}

\title{Thin position and planar surfaces for graphs in the 3-sphere}
\author{Tao Li}
\thanks{Partially supported by an NSF grant}

\address{Department of Mathematics \\
 Boston College \\
 Chestnut Hill, MA 02467}
\email{taoli@bc.edu}
\begin{abstract}
We show that given a trivalent graph in $S^3$, either the graph complement contains an essential almost meridional planar surface or thin position for the graph is also bridge position. This can be viewed as an extension of a theorem of Thompson to graphs. It follows that any graph complement always contains a useful planar surface.   
\end{abstract}
\maketitle

Thin position for a knot is a powerful tool developed by Gabai \cite{Ga}.  Given a knot in thin position, there is a useful planar surface in the knot complement and this planar surface plays an important role in Gabai's proof of property R \cite{Ga} and Gordon-Luecke's solution of the knot complement problem \cite{GL}.  Scharlemann and Thompson later generalized thin position to graphs in $S^3$ and used it in a new proof of Waldhausen's theorem that any Heegaard splitting of $S^3$ is standard, see \cite{ST} and \cite[section 5]{S0}.

In \cite{Th}, Thompson proved that either thin position for a knot is also bridge position or the knot complement contains an essential meridional planar surface.  Wu  improved this result by showing that the thinest level surface is an essential planar surface \cite{Wu}.  In this paper, we generalize Thompson's theorem to graphs in $S^3$ and show that either thin position is also bridge position or there is an essential planar surface in the graph complement and all but at most one of the boundary components of this planar surface are meridians for the graph.  A consequence of this theorem is that there is always a nice planar surface in any graph complement.  The existence of such a nice planar surface is a key in the proof of a theorem in \cite{Li} which says that given a graph $\Gamma$ in $S^3$, if one glues back a handlebody $N(\Gamma)$ to $S^3-N(\Gamma)$ via a sufficiently complicated map, then the resulting closed 3-manifold cannot be $S^3$.

\begin{definition}
Let $N$ be a compact orientable 3-manifold with boundary and let $P$ be an orientable surface properly embedded in $N$.  We say $P$ is \emph{essential} if either $P$ is a compressing disk for $N$ or $P$ is incompressible and $\partial$-incompressible.  We say $P$ is \emph{strongly irreducible} if $P$ is separating, $P$ has compressing disks on both sides, and each compressing disk on one side meets each compressing disk on the other side. $P$ is \emph{$\partial$-strongly irreducible} if 
\begin{enumerate}
	\item every compressing and $\partial$-compressing disk on one side meets every compressing and $\partial$-compressing disk on the other side, and
	\item there is at least one compressing or $\partial$-compressing disk on each side.
\end{enumerate}
Let $\Gamma$ be a graph in $S^3$, $N(\Gamma)$ an open regular neighborhood of $\Gamma$, and $P$ a planar surface properly embedded in  $S^3-N(\Gamma)$.  We say $P$ is \emph{meridional} if every component of $\partial P$ bounds a compressing disk for the handlebody $N(\Gamma)$ and we say $P$ is an \emph{almost meridional planar surface} if all but at most one component of $\partial P$ bounds a compressing disk for $N(\Gamma)$.  Note that a compressing disk for $S^3-N(\Gamma)$ is an almost meridional planar surface by definition.
\end{definition}

Since we are mainly interested in the topology of $N(\Gamma)$ and $S^3-N(\Gamma)$, we may assume our graph $\Gamma$ is trivalent.

\begin{theorem}\label{Tthin}
Let $\Gamma$ be a trivalent graph in $S^3$. Then either
\begin{enumerate}
\item thin position for $\Gamma$ is also bridge position, or
\item $S^3-N(\Gamma)$ contains an essential almost meridional planar surface.
\end{enumerate}
\end{theorem}

The motivation of the paper is the following theorem, which says that a graph complement always contains a nice planar surface.  Such a planar surface turns out to be very useful in \cite{Li}.

\begin{theorem}\label{Tplanar}
Let $\Gamma$ be any graph in $S^3$. Then either 
\begin{enumerate}
\item $S^3-N(\Gamma)$ contains a meridional planar surface that is strongly irreducible and $\partial$-strongly irreducible, or
\item $S^3-N(\Gamma)$ contains an essential almost meridional planar surface, or
\item $S^3-N(\Gamma)$ contains a nonseparating incompressible almost meridional planar surface.  
\end{enumerate}
\end{theorem}
Note that in part (3) of Theorem~\ref{Tplanar}, after some $\partial$-compressions, a nonseparating incompressible planar surface becomes an essential planar surface. Furthermore, it follows from the proof of Theorem~\ref{Tplanar} that the nonseparating planar surface in part (3) of Theorem~\ref{Tplanar} can be chosen to be a punctured disk bounded by an unknotted loop in $\Gamma$ after some edge slides on $\Gamma$

\begin{notation}
Throughout this paper, we denote the interior of $X$ by $\Int(X)$, the closure of $X$ (under the path-metric) by $\overline{X}$, and the number of components of $X$ by $|X|$ for any space $X$.
\end{notation}

The proof of the two theorems is a quick application of the techniques in \cite[section 5]{S0}.  Next we briefly recall some definitions and results in \cite[section 5]{S0}.

\begin{figure} 
\begin{center}
\psfrag{dots}{$\dots\dots$}
\psfrag{du}{$\Delta_u$}
\psfrag{lo}{$\Delta_l$}
\psfrag{D}{$D$}
\psfrag{Y}{$Y$-vertex}
\psfrag{lambda}{$\lambda$-vertex}
\psfrag{(a)}{(a)}
\psfrag{(b)}{(b)}
\psfrag{(c)}{(c)}
\psfrag{(d)}{(d)}
\includegraphics[width=3.5in]{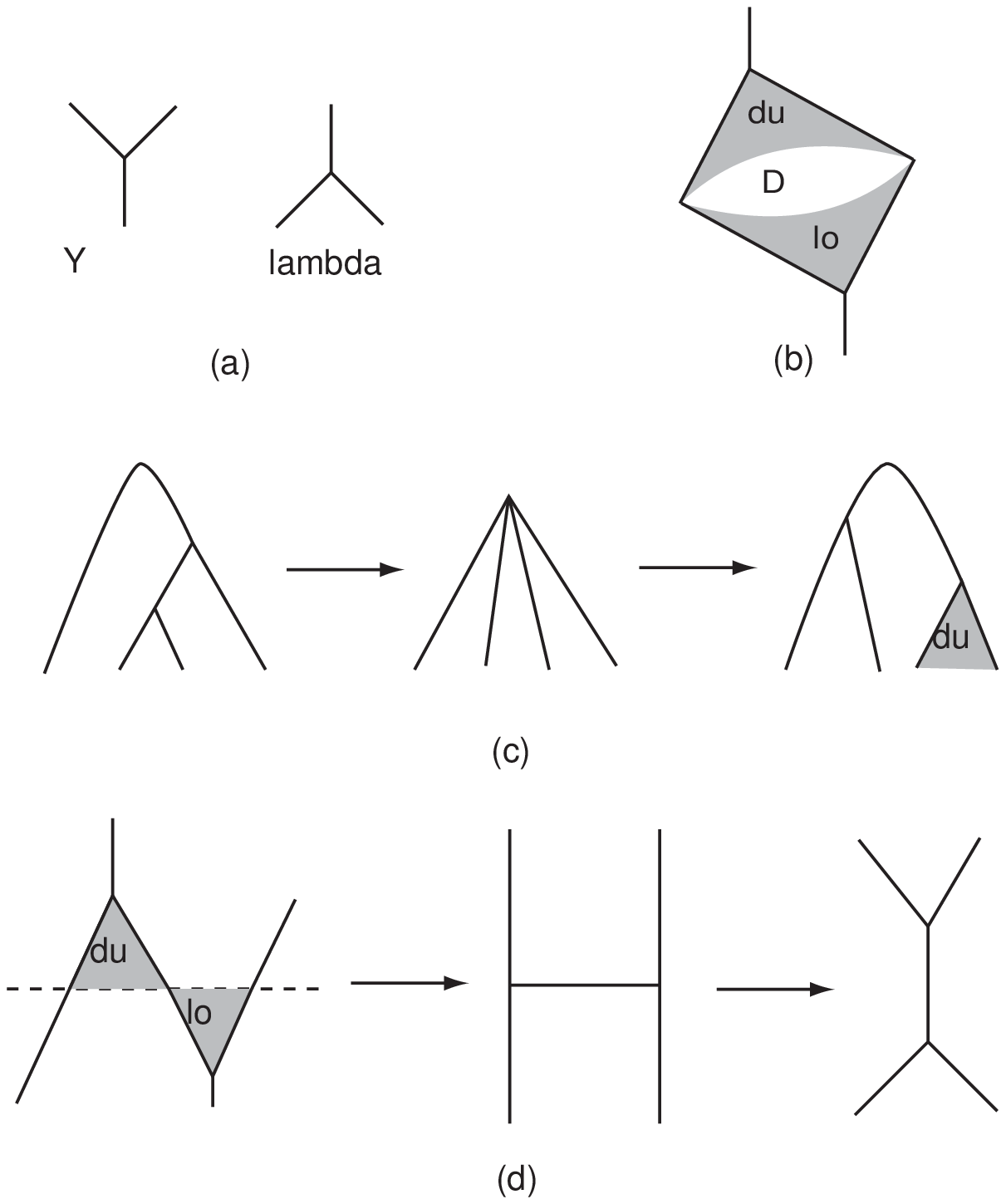}
\caption{}\label{fig}
\end{center}
\end{figure}

We may suppose the trivalent graph $\Gamma$ lies in $S^3-\{\text{two points}\}$ and let $h: S^3-\{\text{two points}\}\to\mathbb{R}$ be the standard height function with each $h^{-1}(p)$  a 2-sphere for any point $p\in\mathbb{R}$.  As in \cite[Definitions 5.1 and 5.2]{S0}, we may assume each vertex of $\Gamma$ is either a $Y$-vertex or a $\lambda$-vertex, as shown in Figure~\ref{fig}(a), and $h|_{edges}$ is a Morse function.  Moreover, the maxima of $\Gamma$ consists of all local maxima of $h|_{edges}$ and all $\lambda$-vertices, and the minima of $\Gamma$ consists of all local minima of $h|_{edges}$ and all $Y$-vertices.  We may also assume the maxima and minima are at different heights and say such a graph $\Gamma$ is in normal form.  $\Gamma$ is in \emph{bridge position} if there is a level sphere that lies above all minima and below all maxima.  The width of $\Gamma$ is defined in \cite[Definition 5.3]{S0} and we say $\Gamma$ is in \emph{thin position} if the height function minimizes the width.  For any $\Gamma$ as above, its width decreases if one pushes a maximum below a minimum, but the width remains the same if one pushes maxima past maxima and minima past minima.  A level sphere $P$ at a generic height is called a \emph{thin sphere} if the adjacent critical height above $P$ is a minimum (possibly a $Y$-vertex) and the adjacent height below $P$ is a maximum (possibly a $\lambda$-vertex).  Clearly  $\Gamma$ is in bridge position if and only if there is no thin sphere.

Next we define upper and lower triples according to \cite[Definition 5.8]{S0}. 
Given $\Gamma$ in normal form and $P$ a level sphere at a generic height, Let $B_u$ and $B_l$ be the two 3-balls above and below $P$ respectively.  Let $E$ be a disk in $S^3-N(\Gamma)$ transverse to $P$ such that $\partial E=\alpha\cup\beta$ where $\alpha$ is an arc properly embedded in $P-N(\Gamma)$, $\partial\alpha=\partial\beta$ and $\beta$ is a nontrivial arc in $\partial\overline{N(\Gamma)}$.  Suppose  the endpoints of $\alpha$ are at different boundary components of $P-N(\Gamma)$.  By shrinking $N(\Gamma)$ back to $\Gamma$, we may view $\alpha$ as an arc in $P$ connecting two punctures of $\Gamma\cap P$.  Let $v$ be one of the points of $\Gamma\cap P$ at an end of $\alpha$ and suppose none of the arc components of $\Int(E)\cap P$ is incident to $v$.  Then $(v,\alpha, E)$ is called an \emph{upper triple} (resp. a \emph{lower triple}) if a small product neighborhood of $\alpha$ in $E$ lies in $B_u$ (resp. $B_l$).  

As in \cite{S0, ST}, if a thin sphere admits an upper or a lower triple $(v,\alpha, E)$, then one can perform an edge slide or a broken edge slide (see \cite[Section 1]{ST} for details) along the disk $E$.  This move however does not reduce the width.  Nonetheless, another combinatorial measure $W_P(\Gamma)$ is introduced in the proof of \cite[Lemma 5.10]{S0} for a thin sphere $P$ and it is shown in \cite{S0} that the (broken) edge slide reduces $W_P(\Gamma)$.  Moreover, the proof of \cite[Lemma 5.10]{S0} shows that if there is always an upper or a lower triple for a thin sphere, then when the process stops, the width is reduced and the thin position for $\Gamma$ is also bridge position.  Although the setting in \cite{S0, ST} is that $N(\Gamma)$ is a handlebody in a Heegaard splitting of $S^3$, the only assumption one needs in the proof of \cite[Lemma 5.10]{S0} is that a thin sphere always admits an upper or a lower triple, see \cite[Lemma 5.9]{S0}.  We now summerize Scharlemann's result as the following lemma.

\begin{lemma}[Scharlemann \cite{S0}]\label{LSch}
Let $\Gamma$ be a trivalent graph in $S^3$ as above.  Then either
\begin{enumerate}
\item thin position for $\Gamma$ is also bridge position, or
\item there is a thin sphere for $\Gamma$ that does not admit any upper or lower triple.
\end{enumerate}
\end{lemma}

Now we are in position to prove Theorem~\ref{Tthin}.

\begin{proof}[Proof of Theorem~\ref{Tthin}]
Since a compressing disk for $S^3-N(\Gamma)$ is an essential almost meridional planar surface, we may suppose $S^3-N(\Gamma)$ has incompressible boundary.  By Lemma~\ref{LSch}, we may assume that there is a thin sphere $P$ that does not admit any upper or lower triple. 

Let $\gamma$ be a boundary component of the planar surface $P-N(\Gamma)$, and let $D$ be an embedded disk in $S^3-N(\Gamma)$ transverse to $P$ with $\partial D=\alpha\cup\beta$, where $\alpha$ is an arc properly embedded in $P-N(\Gamma)$ with both endpoints in $\gamma$, $\beta\subset\partial\overline{N(\Gamma)}$ and $\partial\alpha=\partial\beta$.   Then $\alpha$  cuts $P-N(\Gamma)$ into two planar subsurfaces $P_1$ and $P_2$.  Suppose $\alpha$ is an essential arc in $P-N(\Gamma)$, i.e. neither $P_1$ nor $P_2$ is a disk, and suppose $\beta$ is nontrivial, i.e., $\beta$ is not isotopic (in $\partial\overline{N(\Gamma)}$) relative to $\partial\beta$ to a subarc of $\gamma$.  We say the disk $D$ is a \emph{good disk based at $\gamma$} if for some $i$ ($i=1$ or 2), $\Int(D)\cap P_i$ consists of simple closed curves (if not empty), and we call this $P_i$ the good planar surface associated to $D$.  We call the boundary component of $P_i$ that contains $\alpha$ the outer boundary component of $P_i$ and all other boundary components inner boundary components.  We say the good disk $D$ described above is \emph{innermost with respect to $P_i$} if  there is no good disk $D'$ such that the good planar surface associated to $D'$ is a nontrivial subsurface of $P_i$. 

Suppose $D$ is a good disk based at $\gamma$ and $\Int(D)\cap P_1$ consists of simple closed curves, then one can always find an innermost good disk based at either $\gamma$ or an inner boundary component of $P_1$.  So we may suppose $D$ is innermost with respect to $P_1$.  We claim that $P_1$ is incompressible.  Suppose the claim is not true and there is a nontrivial simple closed curve $\delta$ in $P_1$ such that $\delta$ bounds an embedded disk $\Delta$ in $S^3-N(\Gamma)$.  Since $P_1$ is planar, $\delta$ cuts $P_1$ into two subsurfaces and let $P_\delta$ be the subsurface whose boundary consists of $\delta$ and some inner boundary components of $P_1$. Thus $P_\delta\cup\Delta$ is a meridional planar surface properly embedded in $S^3-N(\Gamma)$.  We may choose $\delta$ to be innermost in the sense that $P_\delta$ is incompressible, i.e. the punctured 2--sphere $S=P_\delta\cup\Delta$ is incompressible in $S^3-N(\Gamma)$. Suppose $S$ is $\partial$-compressible and let $E$ be a $\partial$-compressing disk fro $S$ with $\partial E=\alpha_E\cup\beta_E$, $\alpha_E\subset S$, $\beta_E\subset\partial\overline{N(\Gamma)}$ and $\partial\alpha_E=\partial\beta_E$. Since $S=P_\delta\cup\Delta$ and $\Delta$ is a disk, we may assume $\alpha_E$ is properly embedded in $P_\delta$.  If $\partial\alpha_E$ lies in the same component of $\partial P_\delta-\delta$, since $P_\delta\subset P_1$ and $E\cap P_\delta=\alpha_E$, $E$ is a good disk based at an inner boundary component of $P_1$, which contradicts that $D$ is an innermost good disk with respect to $P_1$.  Thus the two endpoints of $\alpha_E$ lie in different components of $\partial P_\delta-\delta$. By shrinking $N(\Gamma)$ back to $\Gamma$, we may view $\alpha_E$ as an arc connecting two punctures of $\Gamma\cap P$.  Let $v$ be a point of $\Gamma\cap P$ at an end of $\alpha_E$.  Since $E\cap P_\delta=\alpha_E$, $(v,\alpha_E,E)$ is an upper or a lower triple for $P$, which contradicts our earlier assumption.  Thus no such disk $E$ exists and $S=P_\delta\cup\Delta$ is also $\partial$-incompressible.  Hence $S$ is an essential meridional planar surface in $S^3-N(\Gamma)$ and Theorem~\ref{Tthin} holds.  Therefore we may assume $P_1$ is incompressible and this implies that $\Int(D)\cap P_1=\emptyset$ after some isotopy.  So, $P_1\cup D$ is an incompressible planar surface properly embedded in $S^3-N(\Gamma)$.   

Next we show $P_1\cup D$ is also $\partial$-incompressible.  Suppose there is a $\partial$-compressing disk $E'$ for $P_1\cup D$ and suppose $\partial E'=\alpha'\cup\beta'$ where $\alpha'$ is a properly embedded essential arc in $P_1\cup D$, $\beta'\subset\partial\overline{N(\Gamma)}$ and $\partial\alpha'=\partial\beta'$.  Since $D$ is a disk, we can isotope $E'$ so that $\alpha'\subset P_1$.  Similar to the argument above, if the endpoints of $\alpha'$ lie in the same boundary component of $P_1$, then since $E'\cap P_1=\alpha'$, $E'$ is a good disk, which contradicts that the good disk $D$ is innermost with respect to $P_1$.  Suppose the endpoints of $\alpha'$ lie in different boundary components of $P_1$.  Then there is an inner boundary component of $P_1$, say $\gamma'$, containing an endpoint of $\alpha'$.  After shrinking $N(\Gamma)$ to $\Gamma$, $\gamma'$ becomes a point $v'$ in $\Gamma\cap P$ and $\alpha'$ becomes an arc with endpoints in different punctures of $\Gamma\cap P$.  Since $E'\cap P_1=\alpha'$, $(v',\alpha',E')$ is an upper or a lower triple, contradicting our assumption at the beginning of the proof.  This means $P_1\cup D$ must be $\partial$-incompressible.  Thus $P_1\cup D$ is an essential almost meridional planar surface in $S^3-N(\Gamma)$ and Theorem~\ref{Tthin} holds.  Therefore we may assume the thin sphere $P$ does not admit any good disk based at any boundary component of $P-N(\Gamma)$.

Suppose $P-N(\Gamma)$ is compressible in $S^3-N(\Gamma)$ and let $\delta$ be a curve in $P-N(\Gamma)$ that bounds a compressing disk $\Delta$.  Let $U$ be a subdisk of the 2--sphere $P$ bounded by $\delta$. Similar to the argument on $P_1$ above, we may choose $\delta$ to be innermost in the sense that $U-N(\Gamma)$ is incompressible and $Q=\Delta\cup U-N(\Gamma)$ is an incompressible planar surface in $S^3-N(\Gamma)$.  Suppose $Q$ is $\partial$-compressible and let $E''$ be a $\partial$-compressing disk for $Q$ with $\partial E''=\alpha''\cup\beta''$, $\alpha''\subset Q$, $\beta''\subset\partial\overline{N(\Gamma)}$ and $\partial\alpha''=\partial\beta''$.  Similar to the argument above, if $\partial\alpha''$ lies in the same component of $\partial Q$, then $E$ is a good disk based at a boundary component of $P-N(\Gamma)$, which contradicts our assumption above.  If the endpoints of $\alpha''$ lie in different components of $\partial Q$, then $\alpha''$ and $E''$ yield an upper or a lower triple.  Thus $Q$ must be $\partial$-incompressible and hence $Q$ is an essential meridional planar surface and Theorem~\ref{Tthin} holds.
\end{proof}

If thin position for a knot is also bridge position, then it is easy to see that the punctured bridge sphere is strongly irreducible and $\partial$-strongly irreducible in the knot complement.  Theorem~\ref{Tplanar} can be viewed as an extension of this observation to graphs.

\begin{proof}[Proof of Theorem~\ref{Tplanar}]
Since the theorem is about $S^3-N(\Gamma)$, we may assume $\Gamma$ is trivalent.  As in the proof of Theorem~\ref{Tthin}, we may assume $\partial\overline{N(\Gamma)}$ is incompressible in $S^3-N(\Gamma)$.  By Theorem~\ref{Tthin}, either $S^3-N(\Gamma)$ contains a desired essential almost meridional planar surface or thin position for $\Gamma$ is also bridge position.  Suppose thin position for $\Gamma$ is also bridge position and we consider the bridge 2-sphere $S$ for $\Gamma$.  Then $P=S-N(\Gamma)$ is a planar surface in $S^3-N(\Gamma)$ and our goal is to show that either $P$ is strongly irreducible and $\partial$-strongly irreducible or we can construct a nonseparating incompressible almost meridional planar surface in $S^3-N(\Gamma)$.  

The bridge 2-sphere $S$ divides $S^3$ into an upper 3-ball $B_u$ and a lower 3-ball $B_l$.  If $S-\Gamma$ is compressible in $B_u$, then we can maximally compress $S-\Gamma$ in $B_u$ and obtain a collection of mutually disjoint 3-balls $E_1,\dots, E_n$.  Since $S$ is a bridge 2-sphere for $\Gamma$, each $\Gamma_i=\Gamma\cap E_i$ is a connected and unknotted tree.  Each $\Gamma_i$ can be isotoped (relative to the boundary) into a graph in $\partial E_i$.  If one collapse all the edges of $\Gamma_i$ that are not incident to $\partial E_i$ into a single vertex, then $\Gamma_i$ becomes a cone over the points $\Gamma_i\cap\partial E_i$ see \cite[Figure 11]{S0}.  So $\partial E_i-N(\Gamma)$ is parallel to $E_i\cap \partial\overline{N(\Gamma)}$.  Furthermore, $S-\Gamma$ is incompressible in $B_u$ if and only if $\Gamma\cap B_u$ is connected.

If $S-\Gamma$ is incompressible in $B_l$, then $\Gamma\cap B_l$ is a connected unknotted tree.  As in the discussion above, $S-N(\Gamma)$ is parallel to $B_l\cap\partial\overline{N(\Gamma)}$.  Now we consider $B_u$.  By the discussion above, we can find an arc $\alpha$ properly embedded in $\partial E_i-N(\Gamma)$ (where $E_i$ is a 3-ball as above) such that the two endpoints of $\alpha$ lie in different boundary components of $\partial E_i-N(\Gamma)$.  Since $\partial E_i-N(\Gamma)$ is parallel to $E_i\cap\partial\overline{N(\Gamma)}$, $\alpha$ is parallel to an arc $\beta_u$ in $E_i\cap\partial\overline{N(\Gamma)}$.   By the construction of $E_i$, we may view $\alpha$ as an arc in $P$, and there is a disk $D_u$ properly embedded in $B_u-N(\Gamma)$ with $\partial D_u=\alpha\cup\beta_u$, $\partial\alpha=\partial\beta_u$, and $\beta_u\subset B_u\cap\partial\overline{N(\Gamma)}$.  Moreover, since $\Gamma\cap B_l$ is a connected unknotted tree and $\partial B_l-N(\Gamma)$ is parallel to $B_l\cap\partial\overline{N(\Gamma)}$, $\alpha$ is parallel to an arc $\beta_l$ in $B_l\cap\partial\overline{N(\Gamma)}$.  Thus $\beta_u\cup\beta_l$ is a nonseparating simple closed curve in $\partial\overline{N(\Gamma)}$ bounding a disk in $S^3-N(\Gamma)$.  This contradicts our assumption at the beginning that $S^3-N(\Gamma)$ has incompressible boundary.  Thus $S-\Gamma$ and $P$ must be compressible on both sides.

If there are compressing disks $D_u$ and $D_l$ for $S-\Gamma$ in $B_u$ and $B_l$ respectively with $D_u\cap D_l=\emptyset$, then $D_u$ and $D_l$ together with subdisks of $S$ bound two disjoint 3-ball $E_u$ and $E_l$ in $B_u$ and $B_l$ respectively.  Then we can push all the maxima in the $E_u$ into $B_l$ and push all the minima in $E_l$ into $B_u$. Since maxima are pushed below minima, this operation clearly reduces the width and it contradicts that $\Gamma$ is in thin position.  Hence $S-\Gamma$ and $P$ must be strongly irreducible.

The last step is to show that $P$ is also $\partial$-strongly irreducible.  Suppose $P$ is not $\partial$-strongly irreducible, then there is a $\partial$-compressing disk on one side of $P$ disjoint from a compressing or $\partial$-compressing disk on the other side.  Suppose $\Delta_u$ is a $\partial$-compressing disk for $P$ in $B_u$ and $\Delta_l$ is either a compressing or a $\partial$-compressing disk for $P$ in $B_l$ with $\Delta_u\cap\Delta_l=\emptyset$.  Suppose $\partial\Delta_u=\alpha_u\cup\beta_u$ with $\alpha_u\subset P$, $\beta_u\subset B_u\cap\partial\overline{N(\Gamma)}$ and $\partial\alpha_u=\partial\beta_u$.  Next we show that one can choose $\Delta_u$ so that the two endpoints of $\alpha_u$ lie in different components of $\partial P$. 

Suppose $\partial\alpha_u$ lies in the same component $\gamma$ of $\partial P$.  Then $\alpha_u$ must be separating in $P$ and suppose $\alpha_u$ cuts $P$ into two planar subsurfaces $P_1$ and $P_2$.  Since $\Delta_l\cap P$ is either a simple closed curve or a simple arc and $\Delta_u\cap\Delta_l=\emptyset$, we may suppose $\Delta_l\cap P\subset P_2$.  Let $Q$ be the component of $B_u\cap\partial\overline{N(\Gamma)}$ that contains $\gamma$.  Then $\beta_u$ is an arc properly embedded in the planar surface $Q$ with $\partial\beta_u\subset\gamma$. Hence $\beta_u$ is separating in $Q$.  Since $\Gamma\cap B_u$ is unknotted, by the discussion on $\Gamma_i$ and $E_i$ above, at least one component of $\partial Q$, say $\gamma'$, lies in $\Int(P_1)$.  We can find an arc $\beta_u'$ properly embedded in $Q$ with one endpoint in the arc $\gamma\cap\partial P_1$ and the other endpoint in $\gamma'$ and $\beta_u'\cap\beta_u=\emptyset$.  Since $\Gamma\cap B_u$ is unknotted, by the discussion on $\Gamma_i$ and $E_i$ above, $\beta_u'$ is parallel to an arc $\alpha_u'$ properly embedded in $P_1$ with $\partial\alpha_u'=\partial\beta_u'$. Thus $\alpha_u'\cup\beta_u'$ bounds another $\partial$-compressing disk $\Delta_u'$ for $P$ in $B_u$.  The two endpoints of $\alpha_u'$ lie in different components of $\partial P$.  As $\Delta_l\cap P\subset P_2$ and $\alpha_u'\subset P_1$, $\Delta_u'\cap\Delta_l=\emptyset$. Therefore after replacing $\Delta_u$ by $\Delta_u'$ if necessary, we may assume the two endpoints of $\alpha_u$ lie in different components of $\partial P$.  Similarly, if $\Delta_l$ is a $\partial$-compressing disk, we may also assume the endpoints of $\alpha_l=\partial\Delta_l\cap P$ lie in different components of $\partial P$.

Now we shrink $N(\Gamma)$ to $\Gamma$ and to simplify notation, we still use $\Delta_u$ and $\Delta_l$ to denote the disk after the operation.  So $\alpha_u$ is an arcs in $S$ whose endpoints are different punctures of $S\cap\Gamma$.  Since $\Gamma$ is in bridge position, after some edge slides as shown in Figure~\ref{fig}(c) and \cite[Figure 11]{S0}, we may assume $\Delta_u$ is a triangle and $\beta_u$ contains exactly one vertex of $\Gamma$, which is a $\lambda$-vertex.  Similarly, if $\Delta_l$ is also a $\partial$-compressing disk, we may also assume $\Delta_l$ is a triangle whose boundary contains one $Y$-vertex. 

If $\Delta_l$ is a compressing disk, then $\Delta_l$ and a subdisk of $S$ bound a 3-ball $E_l\subset B_l$.  Since $\Delta_u\cap\Delta_l=\emptyset$, we may choose $E_l$ so that $\Delta_u\cap E_l=\emptyset$. We can push $\Delta_u$ and the corresponding $\lambda$-vertex below $S$ and push all the minima in $E_l$ above $S$.  Since a maximum (the $\lambda$-vertex) is pushed below a minimum, this operation reduces the width of $\Gamma$ and it contradicts that $\Gamma$ is in thin position.  So we may suppose $\Delta_l$ is also a $\partial$-compressing disk and assume $\Delta_l$ is triangle as above. Let $\alpha_l=\Delta_l\cap S$. After collapsing $N(\Gamma)$ back to $\Gamma$, $\alpha_u$ and $\alpha_l$ may have common endpoints, but $\Int(\alpha_u)\cap\Int(\alpha_l)=\emptyset$.  We have the following 3 cases to discuss.

The first case is that $\alpha_u\cap\alpha_l=\emptyset$, i.e., $\alpha_u$ and $\alpha_l$ have different endpoints.  Then we can push $\Delta_u$ and the corresponding $\lambda$-vertex below $S$ while pushing $\Delta_l$ and the corresponding $Y$-vertex above $S$.  This reduces the width of $\Gamma$, a contradiction.

The second case is that $\alpha_u$ and $\alpha_l$ share exactly one endpoint.  Then the operation as shown in Figure~\ref{fig}(d) and \cite[Figure 12(a,b)]{S0} reduces the width of $\Gamma$.

The last case is that $\partial\alpha_u=\partial\alpha_l$, see Figure~\ref{fig}(b) and \cite[Figure 12(c)]{S0}.  Since $S^3-N(\Gamma)$ has incompressible boundary, $\alpha_u$ and $\alpha_l$ are not isotopic in $S-\Gamma$.  So $\alpha_u\cup\alpha_v$ bounds a disk $D$ in $S$ and $\Int(D)\cap\Gamma\ne\emptyset$.  We consider the disk $\Delta=\Delta_u\cup D\cup\Delta_l$.  Let $Z$ the set of punctures in $\Int(D)\cap\Gamma$. The punctured disk $\Delta-Z$ corresponds to an almost meridional planar surface $Q$ properly embedded in $S^3-N(\Gamma)$. By our construction, $Q$ is nonseparating.  Next we show that $Q$ is incompressible.  Suppose $Q$ is compressible.  Then there must be a curve $\delta$ in $\Delta-Z$ that bounds a compressing disk for $\Delta-Z$. After isotopy, we may assume $\delta\subset\Int(D)-\Gamma$.  Let $D_\delta$ be the compressing disk bounded by $\delta$ and we may assume $|D_\delta\cap S|$ is minimal among all compressing disks bounded by $\delta$. Let $\delta'$ be a component of $D_\delta\cap S$ that is innermost in $D_\delta$.  Let $D_\delta'$ be the subdisk of $D_\delta$ bounded by $\delta'$.  Then $D_\delta'$ is a compressing disk for $S-\Gamma$ in either $B_u$ or $B_l$. Moreover $D_\delta'$ is disjoint from both $\Delta_u$ and $\Delta_l$.  So similar to the case that $\Delta_l$ is a compressing disk above, we can push $D_\delta'$ and either $\Delta_u$ or $\Delta_l$ to reduce the width of $\Gamma$.  Thus $Q$ must be incompressible in $S^3-N(\Gamma)$ and part (3) of Theorem~\ref{Tplanar} holds.  
\end{proof}

\end{document}